\documentclass[10pt]{amsart}
\usepackage{url,graphicx}
\usepackage[dvips]{epsfig}
\usepackage{latexsym,color}
\usepackage{amscd,amssymb,verbatim}
\newcommand{\be}{\begin{enumerate}}
\newcommand{\ee}{\end{enumerate}}

\newcommand{\aut}{{\text{\rm Aut}\,}}

\newcommand{\cd}{{\mathcal D}}

\newcommand{\cs}{{\mathcal S}}
\newcommand{\rr}{{\mathbb R}}

\newcommand{\da}{\Delta}
\newcommand{\hra}{\hookrightarrow}
\newcommand{\idm}{{\textrm{id}}}
\newcommand{\ra}{\rightarrow}
\newcommand{\nin}{\noindent}

\newcommand{\op}{{\mathfrak o\mathfrak p}}
\newcommand{\pr}{\noindent{\bf Proof. }}
\newcommand{\sm}{\setminus}
\newcommand{\str}{\Sigma_{G,w}}
\newcommand{\cstr}{\overline{\Sigma}_{G,w}}
\newcommand{\bo}{\partial}
\newcommand{\bu}{\bo^\bullet}
\newcommand{\bd}{\bo_\bullet}

\newcommand{\uv}{{\mathfrak{v}}}
\newcommand{\uc}{{\mathfrak c}}
\newcommand{\ucv}{{\mathfrak {cv}}}
\newcommand{\ub}{{\mathfrak b}}

\newcommand{\spr}{\textrm{SP}}
\newcommand{\Br} {\textrm{Br}}
\newcommand{\rp}{{\mathbb R}{\mathbb P}}

\newcommand{\mgn}{MG_n}
\newcommand{\mggn}{MG_{g,n}}
\newcommand{\tmn}{TM_{n}}
\newcommand{\tmgn}{TM_{g,n}}
\newcommand{\mgon}{MG_{1,n}}
\newcommand{\tmon}{TM_{1,n}}
\newcommand{\dz}{{\mathbb Z}}
\newcommand{\zz}{{\mathbb Z}_2}

\newcommand{\hocolim}{\textbf{hocolim}\,}
\newcommand{\hocolimm}{\textbf{\em{hocolim}}\,}
\newcommand{\proj}{\textrm{proj}}
\newcommand{\lra}{\longrightarrow}
\newcommand{\Ree}{\textrm{Re}\,}

\newcommand{\ti}{\tilde}
\newcommand{\wti}{\widetilde}

\newtheorem{thm}{Theorem}[section]
\newtheorem{df}[thm]{Definition}

\newtheorem{crl}[thm]{Corollary}
\newtheorem{prop}[thm]{Proposition}
\newtheorem{conj}[thm]{Conjecture}

\numberwithin{equation}{section}
\numberwithin{figure}{section}

\begin{document}

\title
[Topology of moduli spaces of  tropical curves]
{Topology of moduli spaces of  tropical curves
with marked points}

\author{Dmitry N. Kozlov}
\address{Department of Mathematics, University of Bremen, 28334 Bremen,
Federal Republic of Germany}
\email{dfk@math.uni-bremen.de}
\thanks {This research was supported by University of Bremen, 
as part of AG CALTOP}

\date\today

\begin{abstract}
In this paper we study topology of moduli spaces of tropical curves of
genus~$g$ with $n$ marked points. We view the moduli spaces as being
imbedded in a~larger space, which we call the {\it moduli space of
  metric graphs with $n$ marked points.} We describe the shrinking
bridges strong deformation retraction, which leads to a~substantial
simplification of all these moduli spaces.

In the rest of the paper, that reduction is used to analyze the case
of genus~$1$. The corresponding moduli space is presented as
a~quotient space of a~torus with respect to the conjugation
$\zz$-action; and furthermore, as a~homotopy colimit over a~simple
diagram. The latter allows us to compute all Betti numbers of that
moduli space with coefficients in~$\zz$.
\end{abstract}

\maketitle

\section{Moduli spaces in tropical geometry}

In this paper we study the moduli spaces of  tropical curves
with marked points from the topological point of view. These spaces
were recently introduced by Mikhalkin in \cite{mi,mi2} as important
gadgets in tropical geometry, see also the work of Sturmfels, e.g.,
see \cite{rst} for a~nice introduction. Mikhalkin's investigation
centers on the tropical geometry of these spaces, going in particular
depth in the case of genus~$0$; here we complement his pioneering work
by focusing exclusively on the topological properties.

Accordingly, we define the moduli spaces as imbedded in a~larger space
which we call {\it moduli space of metric graphs with $n$ marked
  points,} in particular, inheriting the natural topology from that
larger space. We then prove that a~simultaneous contraction of all the
bridges is a~strong deformation retraction. The rigorous proof of this
fact is somewhat technical and requires corresponding rigorosity in
the definition of the topology on the set of the isometry classes
of the metric graphs with $n$ marked points. The main technical
problem is to take care of the symmetries arising from the action of
the automorhism group of the graph. All of this is done in Section~3.

Since all the edges in a~tree are bridges, the shrinking bridges
strong deformation retraction contracts the entire moduli space of
 tropical curves of genus $0$ with $n$ marked points to
a~point. Therefore, the corresponding space is not of much interest
from the topological point of view, as far as our current study is
concerned.\footnote{In should be noted that the topology of the link
of that vertex is interesting, and important in the study of
phylogenetic trees, see~\cite{BHV}.}

The first topologically interesting case, that of genus~$1$, is dealt
with in Section~4. In this framework the shrinking bridges strong
deformation retraction simplifies the analysis of the space
dramatically, reducing it to the quotient of the torus by the
conjugation action of~$\zz$.  This can be done because a~bridge-free
graph of genus~$1$ is isomorphic to a~cycle. At present it seems
difficult to describe the homotopy type of the obtained space in
simple terms. However, it is possible to view it as a~homotopy colimit
of a~simple diagram (actually just a~gluing of two identical mapping
cylinders). We use that presentation to compute the Betti numbers of
the space with coefficients in $\zz$, as well as to make a~conjecture
about the homology groups with integer coefficients.

We refer to \cite{CAT} for the concepts and tools of Combinatorial
Algebraic Topology which are used throughout this paper. The few
notions of the category theory that we use here can also be found in
\cite[Chapters 4 and 15]{CAT} or in \cite{mcl}.

\section{Metric graphs}

\subsection{Graphs and graph homomorphisms} $\,$

\nin All the graphs considered in this paper will be finite and
undirected, however loops and multiple edges are allowed, and will, in
fact, be essential for our investigation. Let us now fix our
notations.

\begin{df}\label{df:graph}
A~{\bf graph} $G$ is a~pair of finite sets $V(G)$ and $DE(G)$ equipped
with set maps $\op:DE(G)\ra DE(G)$ and $\bu:DE(G)\ra V(G)$, such that
\begin{itemize}
\item $\op\circ\op=\idm_{DE(G)}$,
\item the map $\op$ has no fixed points.
\end{itemize}
\end{df}

\begin{figure}[hbt]
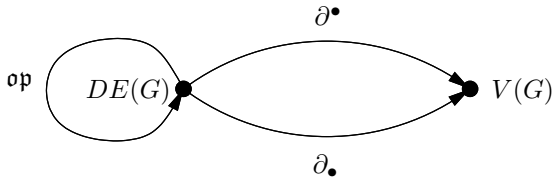

\begin{center}
  \begin{picture}(0,0)%
    \includegraphics{graph1.pstex}%
  \end{picture}%
  \input{graph1.pstex_t}%
  
\end{center}
\caption{A graph $G$ viewed as a~diagram of sets.}
\label{fig:graph}
\end{figure}

When more precise specification is needed, we also write $\op_G$
instead of $\op$.  We shall call $V(G)$ the {\it set of vertices}. We
think of the elements of $DE(G)$ as {\it directed edges}, where the
map $\op$ changes orientation of such a~directed edge to the opposite
one, and the map $\bu$ takes a~directed edge to its source
vertex. Accordingly, we introduce notation $\bd:=\bu\circ\op:DE(G)\ra
V(G)$ for the target vertex map. One can think of a graph as a diagram
of sets over the category with two objects and three non-identity
morphisms, see Figure~\ref{fig:graph}, with the compositions of the
morphisms given by the rules $\op^2=\idm$ and $\bu\circ\op=\bd$.

The condition that $\op$ has no fixed points implies that elements of
$DE(G)$ come in pairs $\{e,\op(e)\}$. These pairs, or equivalence
classes, are the {\it (undirected) edges} of $G$ and we denote the
corresponding set by $E(G):=DE(G)/(e\sim\op(e))$.  Elements $e\in
DE(G)$, such that $\bu e=\bd e$ are called {\it directed loops}, and
the corresponding equivalence classes in $E(G)$ are called {\it
loops}.  For arbitrary vertices $x,y\in V(G)$, we let
\[DE(x,y):=\{e\in DE(G)\,|\,\bu e=x,\,\,\bd e=y\}\] 
denote the set of edges directed from $x$ to $y$. Clearly
$DE(G)=\bigcup_{x,y\in V(G)}DE(x,y)$, and the union is disjoint. The
map $\op$ is a~bijection between $DE(x,y)$ and $DE(y,x)$, and we set
\[E(x,y):=(DE(x,y)\cup DE(y,x))/(e\sim\op(e))\subseteq E(G).\] 
We have $E(x,y)=E(y,x)$, for all $x,y\in V(G)$. In particular
$DE(x,x)$ denotes the set of directed loops at $x$, and
$E(x,x):=DE(x,x)/(e\sim\op(e))$ denotes the set of loops at~$x$.

As an~example, for a~graph $G$ with one vertex and one loop we have
$V(G)=\{v\}$, $DE(v,v)=DE(G)=\{e_1,e_2\}$, with $\op(e_1)=e_2$,
$\op(e_2)=e_1$, $\bu e_1=\bd e_1=\bu e_2=\bd e_2=v$, and
$E(v,v)=E(G)=\{\{e_1,e_2\}\}$, $|E(G)|=1$.

\begin{df} 
For two graphs $G$ and $H$, a~{\bf graph homomorphism} from $G$ to $H$
is simply a~map between corresponding diagrams of sets. In concrete
terms, it consists of two set maps $\varphi_V:V(G)\ra V(H)$ and
$\varphi_{DE}:DE(G)\ra DE(H)$, such that
$\varphi_{DE}\circ\op_G=\op_H\circ\varphi_{DE}$, and
$\varphi_V\circ\bu=\bu\circ\varphi_{DE}$.

A graph homomorphism is called a {\bf graph isomorphism} if the
involved set maps $\varphi_V$ and $\varphi_{DE}$ are bijections.
\end{df}

Since $\varphi_{DE}\circ\op_G=\op_H\circ\varphi_{DE}$ we obtain
induced map $\varphi_E:E(G)\ra E(H)$. For example, the graph with one
vertex and one loop described above has two automorphisms, i.e.,
invertible graph homomorphisms to itself. Both are identity maps on
the sets $V(G)$ and $E(G)$. However, on the set $DE(G)$, one is the
identity map, and the other one swaps the directed edges $e_1$
and~$e_2$.

\subsection{The CW complex $\da(G)$ and the genus of a~graph} $\,$

\begin{figure}[hbt]
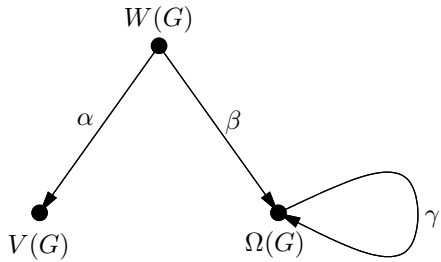

\begin{center}
  \begin{picture}(0,0)%
    \includegraphics{graph2.pstex}%
  \end{picture}%
  \input{graph2.pstex_t}%
  
\end{center}
\caption{The gluing data for the space $\da(G)$.}
\label{fig:wvod}
\end{figure}

\nin We shall now associate a~topological space $\da(G)$ to a~graph
$G$. To avoid making noncanonical choices, and to aid our further
considerations, we would like to think of the space $\da(G)$ as
obtained by gluing together closed intervals corresponding to elements
of $DE(G)$. For this, let $B_e$ denote the closed interval, a~copy of
$[0,1]\subset\rr$, corresponding to the element $e\in DE(G)$, and let
$\Omega(G)$ be the union of all disjoint closed intervals $B_e$.  Let
furthermore $V(G)$ also denote the discrete set of points indexed by
elements of $V(G)$, and let $W(G)$ be the discrete set of points
indexed by the union
\[\{(e,\bu e)\,|\,e\in DE(G)\}\cup\{(e,\bd
e)\,|\,e\in DE(G)\}\subseteq DE(G)\times V(G).\] 
Consider the following maps:
\begin{itemize}
\item a map $\alpha:W(G)\ra V(G)$, defined by $\alpha(e,v):=v$;
\item a map $\beta:W(G)\ra\Omega(G)$, which takes $(e,\bu e)$ to the
  point in $B_e$ corresponding to $0$, and takes $(e,\bd e)$ to the
  point in $B_e$ corresponding to~$1$;
\item a map $\gamma:\Omega(G)\ra\Omega(G)$ which takes a~point in
  $B_e$ corresponding to $x\in[0,1]$ to the point in $B_{\op(e)}$
  corresponding to $1-x$, for all $e\in DE(G)$, and all $x\in[0,1]$.
\end{itemize}
Together with spaces $W(G)$, $V(G)$, and $\Omega(G)$ these maps form
a~diagram shown in Figure~\ref{fig:wvod}, which one can think of as
a~gluing data for $\da(G)$. This is made precise by the following
definition.

\begin{df} \label{df:dag}
For an~arbitrary graph $G$, we let the topological space $\da(G)$ be
the colimit of the diagram shown in Figure~\ref{fig:wvod}.
\end{df} 

We let $q_G:V(G)\cup W(G)\cup\Omega(G)\ra\da(G)$ denote the map
induced by the structural maps from the spaces in the diagram to the
colimit of that diagram. The map $q_G|_{\Omega(G)}$, and hence also
the map $q_G$, is surjective.

Clearly, $\da(G)$ has a~structure of a~$1$-dimensional CW complex,
whose $0$-cells are indexed by the vertices of $G$, $1$-cells are
indexed by the edges of $G$, i.e., by $\op$-invariant pairs of
elements from $DE(G)$, and the attachment maps are given by the
vertex-edge incidences. When appropriate we shall identify vertices of
$G$ with corresponding $0$-cells of $\da(G)$, and edges of $G$ with
corresponding open $1$-cells of $\da(G)$.

A graph homomorphism from $G$ to $H$ induces a~natural diagram map
from the gluing data of $G$ to the gluing data of $H$, and therefore
it also induces a~natural CW map from $\da(G)$ to $\da(H)$, which is
in fact a~homeomorphism when restricted to any open $1$-cell of
$\da(G)$.  We denote both maps by $\varphi_\da$. When $\varphi$ is
a~graph isomorphism, the map $\varphi_\da$ is a~CW isomorphism.

As the last piece of terminology here, the first Betti number of
$\da(G)$ will be called the {\it genus} of $G$, and denoted by
$g(G)$. Clearly $g(G)=|E(G)|-|V(G)|+1$.

\subsection{Metric graphs with marked points.} $\,$ \label{ssect:mgraphs}

\nin To introduce more structure we now vary the lengths of the edges.

\begin{df}  Let $G$ be a~graph. We say that $G$ is a~{\bf metric graph} 
when we are given a~function $l_G:E(G)\ra(0,\infty)$, called the {\bf
edge-length function}. 
\end{df}

The index $G$ will be skipped in $l_G$ whenever it is clear which
graph is considered. We shall also use $l_G$ to denote the corresponding
$\op$-invariant function $l_G:DE(G)\ra(0,\infty)$.

Given a metric graph $(G,l_G)$, there is a standard way to use the
function $l_G$ to turn the topological space $\da(G)$ into a~metric
space, which we now describe. We identify each $B_e$ with the metric
space $[0,l_G(e)]$, instead of the topological space $[0,1]$, with the
standard distance function given by $\ti d(x,y)=|x-y|$, for
$x,y\in[0,l_G(e)]$.

We can now define a~distance function $d$ on $\da(G)$ as follows:
$d(x,x):=0$, for all $x\in\da(G)$, and for $x,y\in\da(G)$, $x\neq y$
we set
\[d(x,y):=\min\sum_{i=1}^n \ti d(x_i,y_i),\]
where the minimum is taken over all $2n$-tuples
$(x_1,y_1,x_2,y_2,\dots,x_n,y_n)$ of points from $\Omega(G)$, such
that points $x_i$ and $y_i$ belong to the same closed interval $B_e$,
for all $i=1,\dots,n$, so the distance $\ti d(x_i,y_i)$ is taken in
this interval, and furthermore $q_G(x_1)=x$, $q_G(y_n)=y$, and
$q_G(y_i)=q_G(x_{i+1})$ for all $i=1,\dots,n-1$.  From now on,
whenever $G$ is a~metric graph, we think of $\da(G)$ as a~metric space
with the standard metric which we just described.

Given a~graph homomorphism from $G$ to $H$, we can adjust the induced
diagram map from the gluing data of $G$ to the gluing data of $H$ to
the metric setting, by taking the map from $B_e$ to $B_{\varphi_E(e)}$
to be the dilation with the scaling factor $l_H(\varphi_E(e))/l_G(e)$.
In the colimit we get the induced map from $\da(G)$ to $\da(H)$, which
we also denote by $\varphi_\da$.

\begin{df}
Let $G$ be a~metric graph, and let $n$ be a~nonnegative integer. We
say that $G$ is a~{\bf metric graph with $n$ marked points}, when we
are given a~function $p_G:[n]\ra \da(G)$ called the {\bf marking
function}.
\end{df}

Here we use the convention $[n]:=\{1,\dots,n\}$ for natural numbers
$n$, and $[0]:=\emptyset$. Formally, metric graphs with $n$ marked
points are given by triples $(G,l_G,p_G)$. Clearly, these generalize
metric graphs, which we can recover by setting $n:=0$.  For
$x\in\da(G)$, we say that {\it $x$ is marked with $p_G^{-1}(x)$}, or
simply that $x$ {\it is marked}, in case that subset of $[n]$ is not
empty. We call a~point $x\in\da(G)$ {\it special} if it is vertex or
a~marked point (or both).

\begin{df}
Two metric graphs $G$ and $H$ with $n$ marked points are said to be
{\bf isometric} if there exists a~graph isomorphism consisting of the
maps $\varphi_V:V(G)\ra V(H)$ and $\varphi_E:DE(G)\ra DE(H)$, such
that we have $l_G=l_H\circ\varphi_E$, and the marked points are mapped
appropriately by the corresponding isometries of the edges, i.e.,
$p_H=\varphi_\da\circ p_G$.
\end{df}

A graph isomorphism $\varphi=(\varphi_V,\varphi_E)$ from $G$ to $H$,
which is also an~isometry of metric graphs, induces an~isometry
$\varphi_\da$ of the corresponding metric spaces.  Isometry of metric
graphs with $n$ marked points is clearly an~equivalence relation. When
$G$ is a metric graph with $n$ marked points we let $[G]$ denote the
corresponding equivalence class (that is the set of all metric graphs
with $n$ marked points which are isometric to $G$). Likewise, for a
set $S$ of metric graphs with $n$ marked points we set
$[S]:=\{[G]\,|\,G\in S\}$.

\subsection{Shrinking edges.} $\,$ \label{ssect:shrink}

\nin Given a~metric graph $G$ with $n$ marked points and $e\in E(G)$,
which is not a~loop, we can define a~new metric graph $H=G/e$ as
follows. Let $\{v,w\}$ be the set of the endpoints of $e$, $v\neq w$,
and let $z$ be a~label which is not in $V(G)$; for reasons which will
become clear shortly, we let the label $z$ be the set $\{v,w\}$
itself. Let furthermore $e_1,e_2\in DE(G)$ be the directed edges
corresponding to $e$. We set $DE(H):=DE(G)\sm\{e_1,e_2\}$ and
$V(H):=(V(G)\sm\{v,w\})\cup\{z\}$. In particular, we see that
$E(H)=E(G)\sm\{e\}$. Let $c:V(G)\ra V(H)$ be the map defined by
$c(v):=c(w):=z$, and $c(x):=x$ for $x\neq v,w$. We now set
$\op_H:=\op_G|_{DE(H)}$ and $\bu_H:=c\circ\bu_G$. In concrete terms,
we have
\[
\begin{array}{rcl}
DE_H(x,y)&:=&DE_G(x,y), \text{ if } z\notin\{x,y\};\\
DE_H(x,z)&:=&DE_G(x,v)\cup DE_G(x,w), \text{ for }x\neq z;\\
DE_H(z,x)&:=&DE_G(v,x)\cup DE_G(w,x), \text{ for }x\neq z;\\
DE_H(z,z)&:=&DE_G(v,v)\cup DE_G(w,w)\cup DE_G(v,w)\cup DE_G(w,v)
\sm\{e_1,e_2\}.
\end{array}
\]
The function $l_H$ is
set to be the restriction of $l_G$ to $E(H)$. Clearly, we have
a~surjective map $\sigma:\da(G)\ra\da(H)$ which ``shrinks'' the edge
$e$, and the marking function $p_H$ is taken to be the
composition~$\sigma\circ p_G$.

More generally, let $S\subseteq E(G)$ be a set of edges which forms
a~subforest of $G$, i.e., the induced graph {\it contains no cycles}
(in particular, the set $S$ contains no loops), and let $DS\subseteq
DE(G)$ be the set of corresponding directed edges. One can then shrink
the set $S$ just like we shrunk a~single edge. More precisely, let
$\Sigma$ be the graph whose set of vertices is $V(G)$ and whose set of
edges is~$S$. The new graph $H=G/S$ is now obtained by taking the
connected components of $\Sigma$ as vertices, and setting
$DE(H):=DE(G)\sm DS$, hence $E(H)=E(G)\sm S$. Let $c:V(G)\ra V(H)$ be
the map taking every vertex of $G$ to the connected component of
$\Sigma$ which contains it. We can then define $\op_H:=\op_G|_{DE(H)}$
and $\bu_H:=c\circ\bu_G$; just like for the case when $S$ consists of
a~single edge. For $C,D\in V(H)$ we now have
\[
\begin{array}{rcl}
DE_H(C,D)&:=&\bigcup_{{x\in V(C)}\atop{y\in V(D)}}DE_G(x,y),\text{ for }C\neq D;\\
DE_H(C,C)&:=&\left(\bigcup_{x,y\in V(C)}DE_G(x,y)\right)\sm DS.
\end{array}
\]
We again let $l_H$ be the restriction of $l_G$ to $E(G)\sm S$, we have
a~surjective shrinking map of topological spaces
$\sigma:\da(G)\ra\da(H)$, and we set $p_H:=\sigma\circ p_G$.

We want to point out a~subtlety related to the edge shrinking. Given
two disjoint sets of edges $S_1$ and $S_2$ we could shrink all these
edges provided that $S_1\cup S_2$ forms a~subforest. If we shrink
first $S_1$ and then $S_2$ we get a~different graph from the one
obtained by shrinking the set $S_1\cup S_2$ right away. This is
because the labels of the vertices will be different. For example, for
a graph with 3 vertices and 2 edges G given by $V(G):=\{a,b\}$,
$E_G(a,b):=\{e\}$, $E_G(b,c):=\{f\}$, shrinking first $e$ and then $f$
yields a~graph with a single vertex labelled $\{\{a,b\},c\}$ and no
edges, while shrinking the entire set $\{e,f\}$ right away yields
a~graph with a single vertex labelled $\{a,b,c\}$ and no
edges. However, it is easy to see that the isometry class $[G/S]$
of the obtained graph does not depend on the order in which we do the
shrinking. We will implicitly use this fact in the future arguments.

\section{The spaces $\mgn$ and $\tmn$, and their deformation retracts}

\subsection{The moduli space of metric graphs with $n$
marked points and its modifications.} $\,$

\nin Let $n$ be a~nonnegative integer, and let $\mgn$ denote the set
of all isometry classes of finite metric graphs with $n$ marked
points. We would like to turn this set into a~topological space. For
this we need to say when two isometry classes of metric graphs with
$n$ marked points ``are close.''

Let $G$ be a metric graph with $n$ marked points. We set $r(G):=\min
d(x,y)$, where the minimum is taken over all pairs of special points
$x,y\in\da(G)$. Note that since the number of special points is
necessarily finite, the minimum is well-defined. We shall refer to the
open interval $(0,r(G)/2)$ as the {\it admissible range} of $G$, this
is the range from which the sizes of the neighborhoods of $[G]$ in
$\mgn$ shall be sampled, and depends only on the isometry class $[G]$,
not on the choice of~$G$.

Let $\varepsilon$ be a~number from the admissible range of $G$. We now
define a set $N_\varepsilon(G)$ as follows: a~metric graph with $n$
marked points $H$ lies in $N_\varepsilon(G)$ if and only if
\begin{enumerate}
\item[(1)] the edges of $H$ of length less than $\varepsilon$ form
  a~subforest;
\item[(2)] the metric graph $G$ can be obtained from $H$ by shrinking
  all the edges of lengths less than $\varepsilon$, as described in
  Subsection~\ref{ssect:shrink}, and by subsequently varying lengths
  of remaining edges and positions of marked points by up
  to~$\varepsilon$.
\end{enumerate}

The latter can be formalized as follows. We say that a metric graph
$G$ with $n$ marked points can be obtained from another metric graph
$H$ with $n$ marked points by varying lengths of edges and the
positions of marked points by up to~$\varepsilon$ if there exists
a graph isomorphism $\varphi$ from $G$ to $H$ such that 
\begin{enumerate}
\item[(1)] for all $e\in E(G)$ we have $|l_G(e)-l_H(\varphi(e))|<\varepsilon$;
\item[(2)] for all $i\in[n]$, we have $d(\varphi_\da(p_G(i)),p_H(i))<\varepsilon$.
\end{enumerate}

Finally, we set
$N_\varepsilon([G]):=[N_\varepsilon(G)]\subseteq\mgn$. Choosing
different representatives of $[G]$ means simply changing labels of
vertices and edges, therefore the set $N_\varepsilon([G])$ is
independent on the particular choice of $G$.

It is easy to see that when $\varepsilon_1$ lies in the admissible
range of $G$, and $\varepsilon_2<\varepsilon_1$, we have
$N_{\varepsilon_2}([G])\subseteq N_{\varepsilon_1}([G])$. This is
because for any $H\in N_{\varepsilon_2}(G)$, and for any edge $e$ of
$H$ we cannot have $\varepsilon_2\leq l_H(e)<\varepsilon_1$, as
otherwise $e$ would correspond to an~edge of $G$, such that
$|l_G(e)-l_H(e)|<\varepsilon_2$, which is impossible since the
inequalities $l_G(e)\geq r(G)$ and $l_H(e)<\varepsilon_1<r(G)/2$ imply
\[l_G(e)-l_H(e)>r(G)-r(G)/2=r(G)/2>\varepsilon_1\geq\varepsilon_2.\]

We are now ready to topologize the set of all isometry classes of
metric graphs with $n$ marked points.

\begin{df}
Let $n$ be a~nonnegative integer. The {\bf moduli space of metric
  graphs with $n$ marked points} is the topological space whose set of
points is given by $\mgn$, and whose topology is generated by the sets
$N_\varepsilon([G])$ as follows: a~subset $X\subseteq\mgn$ is open if
and only if for every $[G]\in X$ there exists $\varepsilon>0$, such
that $N_{\varepsilon}([G])\subseteq X$.
\end{df}

We leave to the reader the verification of the fact that the spaces
$N_{\varepsilon}([G])$ are themselves open. We extend the usage of
$\mgn$ to denote the corresponding topological space as well.

There are various natural modifications of $\mgn$. For example, one
could require the metric graphs to be connected. We denote the
corresponding subspace of $\mgn$ by $\mgn^\uc$. Another, independent
possibility is to require that the marked points are vertices of the
graph. We denote the corresponding subspace of $\mgn$ by $\mgn^\uv$.
Combining, we let $\mgn^\ucv=\mgn^\uc\cap\mgn^\uv$ denote the subspace of
$\mgn$ consisting of the isometry classes of connected metric graphs
with $n$ marks on vertices (as we constructed it, multiple marks are
allowed).

Given a~metric graph $G$ with $n$ marked points, we let $G^\uv$ denote
the metric graph obtained from $G$ by turning all marked points into
vertices (of course, in case they were not vertices already). Clearly,
the isometry class $[G^\uv]$ depends on the isometry class $[G]$ only,
hence the map $\varphi:\mgn\ra\mgn^\uv$ given by $\varphi:[G]\ra[G^\uv]$
is well-defined.

\begin{prop}
The map $\varphi:\mgn\ra\mgn^\uv$ making all marked points into vertices
is a~retraction.
\end{prop}

\pr By construction, we have $\varphi|_{\mgn^\uv}=\idm_{\mgn^\uv}$.
Furthermore, we see that the map $\varphi$ is continuous. Indeed, for
any metric graph $G$ with $n$ marked points the neighborhood
$N_{\varepsilon/2}([G])$ is mapped inside the neighborhood
$N_\varepsilon([G^\uv])$. This is because allowed (for staying in the
neighborhood $N_{\varepsilon/2}([G])$) deformations of the graph $G$:
edge contraction of edges shorter than $\varepsilon/2$, changing the
edge lengths by up to $\varepsilon/2$, shifting a~marked point by at
most $\varepsilon/2$, can all be realized by edge contractions of
edges shorter than $\varepsilon$ and changing edge lengths by at most
$\varepsilon$ in the graph $G^\uv$, when $\varepsilon$ is in the
admissible range of~$G$.  We therefore conclude that $\varphi$ is
a~retraction.
\qed

\subsection{Connected components of $\mgn$.} $\,$

\nin Let us now describe the connected components of $\mgn$. Let $G$ be
a~metric graph with $n$ marked points, and let $G_1,\dots,G_t$ be its
connected components. Assume that $G_i$ has genus $g_i$, for
$i=1,\dots,t$, and let $A_1,\dots,A_t$ be disjoint, possibly empty
sets whose union is $[n]$ (this is like a~set partition, but with
empty sets allowed). The set $\{(g_1,A_1),\dots,(g_t,A_t)\}$ is now
the data which we associate to the graph $G$. It is not difficult to
see that any graph with the same data lies in the connected component
which contains $[G]$. Furthermore, the data of this type, meaning
a~set of tuples $(g_i,A_i)$, such that $g_i\geq 0$ and
$(A_1,\dots,A_t)$ is a~set partition of $[n]$, possibly involving
empty sets, index connected components of $\mgn$.

Consider now a~connected component $C_S$ indexed by the set
\[S=\{(g_1,A_1),\dots,(g_m,A_m),(0,\emptyset),\dots,(0,\emptyset),
(1,\emptyset),\dots,(1,\emptyset),\dots\},\] where the sets $A_i$ are
non-empty. Denote the number of appearances of the tuple
$(k,\emptyset)$ in that set by $n_k$, for $k=0,1,\dots$. By our
assumptions, only finitely many of these are different from~$0$.
Then we have a~homeomorphism 
\[C_S\cong C_{\{(g_1,A_1)\}}\times\dots\times C_{\{(g_1,A_1)\}}\times
\spr^{n_0}(C_{\{(0,\emptyset)\}})\times\dots\times
\spr^{n_i}(C_{\{(i,\emptyset)\}})\times\dots,\] where $\spr^t(X)$
denotes the $t$-fold symmetric product, i.e., the quotient space
$(\underbrace{X\times\dots\times X}_t)/\cs_t$, where the symmetric
group $\cs_t$ acts on the direct product by permutation of its
factors. 

We shall use $\mggn$ to denote the topological space whose points are
the isomorphism classes of connected graphs of genus $g$ with $n$
marked points, which is the same as the connected component
$C_{(g,[n])}$.

\subsection{A stratification of $\mgn$} $\,$

\nin The moduli space $\mgn$ has a~natural stratification. To produce
a~stratum, fix a~graph $G$, and for each $i\in[n]$, fix $w_i$, which
is either a~vertex or an edge of $G$. Now, consider the set of all
isometry classes of metric graphs with $n$ marked points, which have
a~representative $(H,l_H,p_H)$, such that there exists a~graph
isomorphism $\varphi$ from $G$ to $H$, for which the point $p_H(i)$ is
equal to $\varphi(w_i)\in\da(H)$, if $w_i$ is a~vertex, or belongs the
open edge $\varphi_{x,y}(w_i)\subset\da(H)$, if $w_i$ is an~edge with
endpoints $x$ and $y$, for all $i\in[n]$.

This stratum is indexed by the graph $G$ together with the $n$-tuple
$w=(w_1,\dots,w_n)$; we denote it by $\str$. We shall call this
stratification {\it standard}, and we shall its strata the {\it
standard strata}. The stratum does not change if we replace $G$ with
an isomorphic graph, and change the $n$-tuple $w$ accordingly. We
shall implicitly use this fact in our discussion.

We let $\cstr$ denote the closure of the stratum $\str$, and we let
$\bo\str$ denote $\cstr\sm\str$, which we shall call the {\it
boundary} of the stratum. One can see that the boundary of
an~arbitrary stratum $\str$ is a~union of other strata. The indexing
data of these strata can be obtained from $(G,w)$ by a~combination of
the steps of the following two kinds
\begin{enumerate}
\item[(1)] replacing an~edge by one of its endpoints in the
  $n$-tuple~$w$;
\item[(2)] shrinking a~non-loop edge in $G$ and replacing this edge and its
endpoints by the label of the thus obtained vertex in the $n$-tuple~$w$.
\end{enumerate}

In general the strata do not have to be manifolds. Consider, for
example, the stratum $\str$, where $G$ is a~graph with one vertex $v$
and one edge $e$ (which hence must be a~loop), $n=1$, and $w_1=e$. We
can vary the length of the edge in the open interval $(0,\infty)$, and
we can slide the marked point along the edge. Thus potentially the
point moves in the interval $(0,l_G(e))$, however, because of the
symmetry which flips the loop, we need to identify coordinates $x$ and
$l_G(e)-x$, so we can choose a~representative from the half-closed
interval $(0,l_G(x)/2]$. These considerations show that the stratum
$\str$ is homeomorphic to the space $\{(l,x)\in\rr^2\,|\,l>0,\,\,l/2\geq
x>0\}$, which in turn is homeomorphic to the space
$\{(x,y)\in\rr^2\,|\,x>0,\,\,y\geq 0\}$.

On the positive side, as easily seen, the generic points of every
stratum $\str$ form an~open manifold whose dimension is equal to the
number of edges of $G$ plus the number of labels which are edges in
the $n$-tuple~$w$. It follows, that there are infinitely many strata
of dimension $0$; these are indexed by graphs $G$ with no edges, whose
vertices are labeled by disjoint subsets of $[n]$, so that the union
of all labels is $[n]$, which is the same as to index them by sets
$\{A_1,\dots,A_t\}$, where we might have $A_i=\emptyset$, such that
$[n]=\cup_{i=1}^t A_i$, and the union is disjoint.

In general, the space $\mgn$ is somewhat technical to handle directly:
it is infinite dimensional, and an arbitrarily small neighborhood of
each point intersects infinitely many strata; for example, in the case
with no marked points, an arbitrarily small neighborhood of the graph
with one vertex and no edges intersects all strata indexed by trees.
We shall therefore start by performing the shrinking bridges strong
deformation retraction, in order to replace $\mgn$ by a~more
manageable space.

\subsection{The shrinking bridges strong deformation retraction} $\,$

\nin Recall, that an edge $e$ of a graph $G$ is called a~{\it bridge}
if deleting it from the graph $G$ increases the number of connected
components. Equivalently, $e$ is a~bridge if the endpoints of $e$
belong to different connected components of~$G-e$
(cf.\ \cite[p.\ 11]{Die}). Clearly, shrinking a~bridge will neither
change the number of connected components of $G$, nor will it change
the genuses of these connected components. For example, every edge of
a~forest is a~bridge. We shall call a~graph which has no bridges {\it
  bridge-free}, and we shall denote the set of bridges of $G$ by
$\Br(G)$.

Let us now define a~homotopy $\beta:\mgn\times[0,1]\ra\mgn$. Let $G$
be a~metric graph with $n$ marked points, and let
$t\in(0,1]$. Let~$\beta(G,t)$ denote the metric graph with $n$ marked
  points obtained from $G$ by scaling down all the edges in $\Br(G)$
  by the factor~$t$, and adjusting the marking function
  accordingly. For $t=0$ we set $\beta(G,0)$ to be the graph obtained
  from $G$ by shrinking all the bridges, this is allowed since the set
  of all bridges forms a~subforest of~$G$. Up to isomorphism, the
  metric graph $\beta(G,t)$ with $n$ marked points is uniquely
  determined by the isomorphism class of $G$, and by the
  parameter~$t$, hence the assignment $\beta([G],t):=[\beta(G,t)]$ is
  well-defined.

Let $\mgn^\ub$ denote the subspace of $\mgn$ consisting of all the
isometry classes of metric bridge-free graphs $G$ with $n$ marked
points. Clearly, the space $\mgn^\ub$ is a~union of standard strata, and
the map $\beta(-,0)$ surjectivity takes $\mgn$ to $\mgn^\ub$.

\begin{thm} \label{thm:bbh}
The space $\mgn^\ub$ is a~strong deformation retract of the space
$\mgn$. The map $\beta:\mgn\times[0,1]\ra\mgn$ provides
a~corresponding strong deformation retraction.
\end{thm} 

\pr As already mentioned, we have $\beta([G],0)\in\mgn^\ub$, for any
metric graph $G$ with $n$ marked points. Furthermore, it follows
directly from our definition of the map $\beta$, that
$\beta([G],t)=[G]$, for all $[G]\in\mgn^\ub$, and that
$\beta(-,1)=\idm_{\mgn}$. Therefore, to prove that $\beta$ is
an~appropriate strong deformation retraction, it is enough to show
that it is continuous. 

We shall now provide a~direct, albeit somewhat tedious verification.
We take a~point $x\in\mgn$ and a~point $y\in\mgn\times[0,1]$, such
that $\beta(y)=x$, and then show that for every sufficiently small
$\varepsilon$ (how small it needs to be shall depend on $x$ and $y$),
there exists a~neighborhood $N$ of $y$, which is taken by $\beta$
inside of $N_\varepsilon(x)$: $\beta(N)\subseteq N_\varepsilon(x)$.
We use the notations $y=([G],t)\in\mgn\times[0,1]$, where $G$ is
a~metric graph with $n$ marked points, and accordingly
$x=\beta([G],t)=[\beta(G,t)]$. It is technically easier to divide the
argument into considering two separate cases.

\vskip5pt

\nin {\bf Case 1.} {\it We assume that $t>0$.}

\nin This is the easier one of the two cases.  For sufficiently small
$\varepsilon>0$, we look for $\delta_1,\delta_2>0$, such that $\beta$
maps the $y$-neighborhood
$N_{\delta_1}([G])\times[t-\delta_2,t+\delta_2]$ inside of
$N_\varepsilon(x)$. Since $N_{\delta_1}([G])=[N_{\delta_1}(G)]$, and
$N_\varepsilon(x)=N_\varepsilon(\beta([G],t))=
N_\varepsilon([\beta(G,t)])= [N_\varepsilon(\beta(G,t))]$, it is
enough to find $\varepsilon$, $\delta_1$, and $\delta_2$, such that
$\beta$ maps $N_{\delta_1}(G)\times[t-\delta_2,t+\delta_2]$ inside of
$N_\varepsilon(\beta(G,t))$, as long as the conditions on
$\varepsilon$, $\delta_1$, and $\delta_2$ depend only on the
isomorphism class of $G$, not on the specific representative.  In any
case, we assume that $\varepsilon$ is sampled from the admissible
range of~$G$ (which only depends on $[G]$).

Let now $H$ be a~metric graph with $n$ marked points in
$N_{\delta_1}(G)$, and let $\Sigma$ denote the set of edges of $H$ of
length less than $\delta_1$. By our construction, these must form
a~subforest. Let us fix some $\ti t$ in the interval
$[t-\delta_2,t+\delta_2]$. The graph $\beta(H,\ti t)$ is obtained from
$H$ by shrinking all the bridges by a~factor~$\ti t\leq 1$, $\ti t\neq
0$. Therefore, requesting that $\delta_1\leq\varepsilon$ will ensure
that the images of the edges from $\Sigma$ will have length less than
$\varepsilon$.  On the other hand, we want the images of the edges
from $E(H)\sm\Sigma$, i.e., the original edges of $G$, to have lengths
larger than~$\varepsilon$. This can be ensured by requesting that
\[\ti t\cdot\min_{e\in E(G)}l_G(e)>\varepsilon.\]
Therefore, the graph $\beta(G,\ti t)$ can be obtained from the graph
$\beta(H,\ti t)$ by shrinking all the edges of length less
than~$\varepsilon$, and then varying the lengths of the remaining
edges, as well as positions of marked points, by up to
an~$\varepsilon$. The latter follows from the fact that
$\delta_1\leq\varepsilon$, and shrinking some of the edges only
decreases the edge length variation.

Passing on to the whole interval $[t-\delta_2,t+\delta_2]$, we choose
$\varepsilon$ so that 
\[t\cdot\min_{e\in E(G)}l_G(e)/2>\varepsilon,\]
and then take $\delta_1\leq\varepsilon$ and $\delta_2\leq t/2$. This
choice of parameters verifies the continuity of the map $\beta$ at the
point~$y$.

\vskip5pt

\nin {\bf Case 2.} {\it We assume that $t=0$.}

\nin We have $x=\beta([G],0)\in\mgn^\ub$. This time, for sufficiently
small $\varepsilon>0$, we need to find $\delta_1,\delta_2>0$, such
that $\beta$ maps the $y$-neighborhood
$N_{\delta_1}([G])\times[0,\delta_2]$ inside of
$N_\varepsilon(x)$. Just like in the first case, we can drop the
isomorphism brackets, and search for $\varepsilon$, $\delta_1$, and
$\delta_2$, such that $\beta$ maps
$N_{\delta_1}(G)\times[0,\delta_2]$ inside of $N_\varepsilon(\beta(G,0))$.

Let again $H$ be a~metric graph with $n$ marked points in
$N_{\delta_1}(G)$, and let $\Sigma$ denote the set of edges of $H$ of
length less than $\delta_1$. It important to note, that the set of
edges $\Sigma\cup \Br(G)$ forms a~subforest. This is because $\Sigma$
is a~forest, and adding bridges to any forest will not create cycles,
since bridges cannot be a~part of any cycle. Also, as noted before,
we have $\Br(G)\subseteq \Br(H)\subseteq \Br(G)\cup\Sigma$.

The graph $\beta(H,0)$ is obtained from $H$ by shrinking all the
bridges.  This means shrinking all the bridges of $G$, and possibly
some of the edges from~$\Sigma$. Choosing $\varepsilon$ smaller than
$\min_{e\in E(G)}l_G(e)$, and then choosing $\delta_1<\varepsilon$
ensures that the edges of $\beta(H,0)$ whose lengths are less than
$\varepsilon$ are precisely the non-bridges from the set
$\Sigma$. This means that the graph $\beta(G,0)$ can be obtained from
the graph $\beta(H,0)$ by shrinking the edges of length less than
$\varepsilon$, and varying the lengths of other edges, as well as positions
of the marked points, by up to an~$\varepsilon$.

Assume now $\delta_2$ is chosen so that 
\[\delta_2\cdot\max_{e\in E(G)}l_G(e)<\varepsilon,\]
and choose $0<t\leq\delta_2$. The graph $\beta(H,t)$ is obtained from
$H$ by scaling all the bridges down by the factor $t$. The way
$\delta_2$ is chosen, this will scale down all the bridges of $G$, so
that they become shorter than $\varepsilon$. Furthermore, the
conditions $\delta_1<\varepsilon<\min_{e\in E(G)}l_G(e)$ imply that
also all the edges from $\Sigma$ will be shorter than~$\varepsilon$
(since they were shorter than $\varepsilon$ to start with, and some
were additionally shrunk), and that all the non-bridges of $G$ will
not be shorter than~$\varepsilon$. Hence again the graph $\beta(G,0)$
can be obtained from the graph $\beta(H,t)$ by shrinking the edges of
length less than $\varepsilon$, and varying the lengths of other
edges, as well as positions of the marked points, by up to
an~$\varepsilon$.

This finishes the verification of the fact that the homotopy $\beta$
is continuous at the point $y$ in this case.
\qed

\vskip5pt

\subsection{Moduli space of tropical curves of genus $g$ with $n$ marked
points.} $\,$

\nin We now define a subspace of $\mgn^\ucv$ which is of special interest in
tropical geometry and has been the starting point of the current
investigation.

\begin{df}
Let $n$ be a~nonnegative integer, and let $d$ be a~positive number. We
define $\tmn(d)$ to be the subspace of $\mgn^\ucv$ consisting of the
isomorphism classes of all metric graphs $G$ with $n$ marked points,
such that
\begin{enumerate}
\item[(1)] $G$ has no vertices of valency 2;
\item[(2)] $G$ has exactly $n$ leaves\footnote{Generalizing the
  terminology customary for trees, we use the word {\it leaves} to
  denote any vertex of valency $1$, cf.\ \cite[p.\ 13]{Die}.}, and
  these are marked $1$ through $n$;
\item[(3)] the lengths of the edges leading to leaves are equal to $d$. 
\end{enumerate}
\end{df}

Simultaneous dilation of the edges leading to leaves gives
a~homeomorphism between spaces $\tmn(d_1)$ and $\tmn(d_2)$, for
arbitrary positive $d_1$ and $d_2$. Letting $d$ go to $0$ we obtain
yet another homeomorphic space, which it would be natural to denote by
$\tmn(0)$, but for simplicity we just call it $\tmn$. This space can
also be described directly, as is done in the next definition.

\begin{df} \label{df:tmn}
Let $n$ be a~nonnegative integer. We define $\tmn$ to be the subspace
of $\mgn^\ucv$ consisting of the isomorphism classes of all metric
graphs $G$ with $n$ marked points, such that for every vertex of $G$
the sum of its valency with the number of times it is marked should be
at least~3.
\end{df}
\nin The condition in Definition~\ref{df:tmn} just means that every
vertex of valency 2 should be marked, and that every leaf should be
marked at least twice.

It is not a~difficult exercise to see that every two points of $\tmn$
corresponding to metric graphs of the {\it same} genus, can be
connected by a~path inside $\tmn$, whereas obviously, every two points
of $\tmn$ corresponding to metric graphs of {\it different} genus,
cannot be connected by such a~path, not even inside of $\mgn$. Hence
the connected components of $\tmn$ are indexed by nonnegative integers
$g$, corresponding to the genuses of the involved graphs, and we call
them $\tmgn$.\footnote{It has been recently proved that the tropical
  moduli space $\tmon$ is a~strong deformation retract of the graph
  moduli space $\mgon^\uv:=\mgon\cap\mgn^\uv$, see~\cite{tms2}.}

It is important to not that, unlike the spaces $\tmn(d)$, the space
$\tmn$ is a union of the standard strata of $\mgn$.

\begin{crl} \label{crl:tmnb}
Let $\tmn^\ub:=\tmn\cap\mgn^\ub$. Then the space $\tmn^\ub$ is a~strong
deformation retract of~$\tmn$.
\end{crl}
\pr
The space $\tmn$ is obviously closed under the shrinking bridges strong
deformation retraction. Therefore, the Theorem~\ref{thm:bbh} implies
that $\tmn$ strongly deformation retracts to $\tmn^\ub$.
\qed 
\vskip5pt

The space $\tmn^\ub$ can be also described directly: it is the subspace
of $\mgn^\ucv$ consisting of all isomorphism classes of metric graphs
$G$ with $n$ marked points, such that
\begin{enumerate}
\item[(1)] every vertex of $G$ of valency 2 is marked;
\item[(2)] the graph $G$ has no bridges.
\end{enumerate}

The connected components of $\tmn^\ub$ are again indexed by genuses of
the constituting graphs, and we shall use the notation
$\tmgn^\ub:=\tmgn\cap\mgn^\ub$.


\section{The moduli space of the tropical curves of genus $1$}

\subsection{Presentation as a~quotient space} $\,$

\nin Since all the edges of a~tree are bridges, we see that $MG_{0,n}^\ub$
and $TM_{0,n}^\ub$ are just single points. Theorem~\ref{thm:bbh} and
Corollary~\ref{crl:tmnb} imply that the spaces $MG_{0,n}$ and
$TM_{0,n}$ are contractible. In this section we shall focus on the
next interesting case: namely the spaces of connected metric graphs
of genus~$1$ with $n$~marked points.

Let us start by analysing the space $\wti X_n=TM_{1,n}^\ub$. The
bridge-free graphs of genus $1$ are simply cycles, hence, since all
their vertices have valency equal to $2$, they should all be marked,
and $n$ should be at least~$1$. Reversely, all the marked points are
vertices. We can therefore forget about the vertices and just record
the marked points.  

The points of the space $\wti X_n$ can thus be indexed by $n$-tuples
of points on a~circle, whose radius is an arbitrary positive real
number, divided by the action of the orthogonal group on the circle.
Factoring out the radius length of the circle, we have a~homeomorphism
\[\wti X_n\cong X_n\times (0,\infty),\]
for $X_n:=(\underbrace{S^1\times\dots\times S^1}_n)/O(2)$, where $S^1$
denotes the unit circle, and the division is done with respect to the
diagonal action of the orthogonal group $O(2)$, which acts on each
term in the natural way.

Note, that the described action of $O(2)$ on the direct product
$\underbrace{S^1\times\dots\times S^1}_n$ is transitive on the last
coordinate. Let us fix the last coordinate to be $(1,0)\in S^1$. Then
the space $X_n$ can be rewritten as
\begin{equation} \label{eq:xn}
X_n\cong \underbrace{S^1\times\dots\times S^1}_{n-1}/\zz,
\end{equation}
where $\zz$ is the subgroup of $O(2)$ which fixes the point $(1,0)$.
Clearly, the group $\zz$ consists of two elements, and the
non-identity element is an~involution which acts diagonally on the
direct product of $n-1$ copies of $S^1$ by a~reflection about the
$X$-axis. One way to think of this action is to view the points of
$S^1$ as complex numbers with absolute values equal to $1$, in which
case the action is just a~simultaneous conjugation:
\[X_n\cong\{(z_1,\dots,z_{n-1})\,|\,|z_i|=1,\textrm{ for }
i=1,\dots,n-1\}/(z_1,\dots,z_{n-1})\sim(\bar z_1,\dots,\bar z_{n-1}).\]

\subsection{A cubical structure on $X_n$} $\,$ \label{ssect:4.2}

\begin{figure}[hbt]
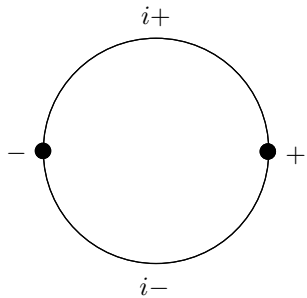

\begin{center}
  \begin{picture}(0,0)%
    \includegraphics{circle.pstex}%
  \end{picture}%
  \input{circle.pstex_t}%
  
\end{center}
\caption{The cell encoding in Subsection~\ref{ssect:4.2}.}
\label{fig:circle}
\end{figure}

\nin Consider the CW structure on a unit circle $S^1$ consisting of
two $0$-cells $(1,0)$ and $(-1,0)$, and two $1$-cells corresponding to
upper and lower semicircles. This induces a~CW structure on the direct
product of $n$ copies of $S^1$ (which is an $n$-torus), where the
cells are simply direct products of cells of the factors. These are
indexed by the ordered $n$-tuples of the cells of the factors and we
shall use the following encoding: we write ``$+$'' for the $0$-cell
$(1,0)$, ``$-$'' for the $0$-cell $(-1,0)$, ``$i+$'' for the upper
semicircle, and ``$i-$'' for the lower semicircle, as shown on
Figure~\ref{fig:circle}. So, for example, $(i-,-,-,i+)$ would denote
a~$2$-cell in the $4$-torus. The open cells are actually cubes, so we
can think of this decomposition as some sort of a~cubical structure on
$\underbrace{S^1\times\dots\times S^1}_n$. In this encoding, the
boundary of a~cell is obtained by replacing symbols $i+$ and $i-$ with
symbols $+$ and~$-$. The number of $d$-cells is ${{n}\choose{d}}2^n$,
as we can put a~plus or a~minus in every coordinate, and then
distribute $d$ symbols ``$i$'' arbitrarily.

Clearly, this cubical structure on $\underbrace{S^1\times\dots\times
  S^1}_{n-1}$ is invariant under the $\zz$-action of simultaneous
conjugation. In our symbolic notations, the action changes the signs
assigned to $i$'s and fixes the other coordinates. In particular, all
$0$-cells are fixed, and all higher dimensional cells come in pairs,
the cells in each pair are being swapped by the involution. This means
that in the quotient $X_n$ we get the induced cubical structure,
consisting of the orbits of the cells of the $(n-1)$-torus. These are
also indexed by the $n$-tuples of the symbols from the set
$\{+,-,i+,i-\}$, with an additional constraint that the first symbol
$i+$ comes before the first symbol $i-$ (if at all).  Hence the number
of vertices in this cubical structure on $X_n$ is $2^{n-1}$, whereas
the number of $d$-cubes, for $n-1\geq d\geq 1$, is
${{n-1}\choose{d}}2^{n-2}$.

\subsection{The small values of $n$} $\,$

\nin The presentation~\eqref{eq:xn} can be used to analyze what
happens for the small values of $n$. The space $X_1$ is just a point:
the empty direct product is to be interpreted as a point here. The
space $X_2\cong S^1/\zz$ is homeomorphic to a~closed interval, e.g.,
$[-1,1]$ if we take the projection of $S^1$ onto the $X$-axis.

Let us consider $X_3\cong(S^1\times S^1)/\zz$. We have a~described
a~cubical structure on this space which has $4$ vertices, $4$ edges
and $2$ squares. The two squares are indexed with $(i+,i+)$ and
$(i+,i-)$ and it is easy to see that $X_3$ is obtained by gluing these
$2$ squares together along their entire boundaries. Hence we conclude
that $X_3$ is homeomorphic to a~$2$-dimensional sphere.

\begin{figure}[hbt]
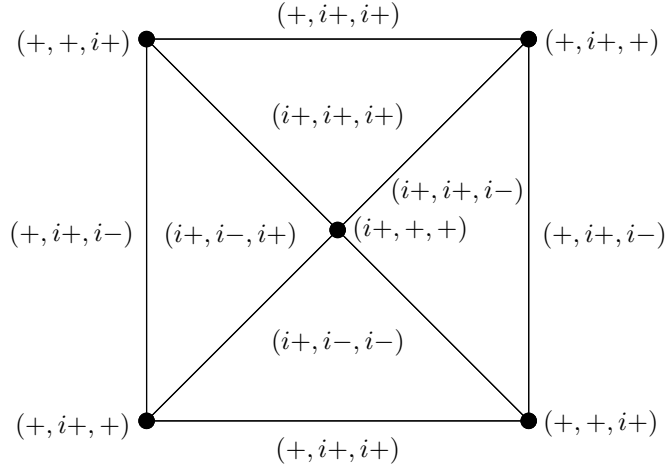

\begin{center}
  \begin{picture}(0,0)%
    \includegraphics{link1.pstex}%
  \end{picture}%
  \input{link1.pstex_t}%
  
\end{center}
\caption{The link of the vertex $(+,+,+)$ in $X_4$.}
\label{fig:link}
\end{figure}

As the last case, let us consider $X_4\cong(S^1\times S^1\times
S^1)/\zz$.  This is a~space glued together from $4$ cubes; it has $8$
vertices, $12$ edges and $12$ squares. All of the cubes share the same
set of vertices, and it is easy to see that all points, except for
these $8$ vertices, have neighborhoods which are homeomorphic to open
balls in $\rr^3$. However, the space $X_4$ is not a~manifold. To see
this, let us take a~look at the neighborhoods of these
vertices. Clearly, it does not matter which one we take, so let us
take the vertex $v=(+,+,+)$. As an open neighborhood of $v$ we can
take a cone with an apex in $v$ over the link of $v$. The link of
a~vertex in a~cubical complex is always a~simplicial complex. The one
we have here has $3$ vertices, $6$ edges, and $4$
triangles. A~presentation of this simplicial complex is given on
Figure~\ref{fig:link}, from which it is clear that this link is
homeomorphic to $\rp^2$. This shows that $X_4$ fails to be a~manifold
at these points.

\subsection{Surgery presentation} $\,$

\nin Let us now generalize our description of $X_4$ to the general
case. In particular, we shall see that for $n\geq 4$ the vertices in
the cubical structure on $X_n$ are singularities, and if they are
removed, we are left with an~$(n-1)$-dimensional manifold.

We start with by viewing the $n$-torus $T^n$ in the standard way as
the quotient space of $\rr^n$ divided by the group action of
$\dz^n=\langle g_1,\dots,g_n\rangle$, given by 
\[g_i:(x_1,\dots,x_{i-1},x_i,x_{i+1},\dots,x_n)\mapsto
(x_1,\dots,x_{i-1},x_i+2,x_{i+1},\dots,x_n).\] Since $X_{n+1}\cong
T^n/\zz$, as described above, we have $X_n\cong\rr^n/G_n$, where $G_n$ 
is a group defined by 
\begin{equation} \label{eq:gn}
G_n=\langle\gamma,g_1,\dots,g_n\,|\,g_i\circ\gamma=\gamma\circ g_i^{-1},\,\,
g_i\circ g_j=g_j\circ g_i,\text{ for }i,j\in [n]\rangle,
\end{equation}
with the action of $\gamma$ on $\rr^n$ given by 
\[\gamma:(x_1,\dots,x_i,\dots,x_n)\mapsto
(-x_1,\dots,-x_i,\dots,-x_n).\] 

It follows from the presentation~\eqref{eq:gn} that the group $G_n$ is
the semidirect product $\dz^n\times_\phi\zz$, with the group
homomorphism $\phi:\zz\ra\aut(\dz^n)$ given by
$\phi(\gamma)(g)=g^{-1}$, where $\gamma$ is the nontrivial element of
$\zz$ and $g\in\dz^n$ is arbitrary.  In particular, the elements of
$G_n$ can be uniquely presented either as $\gamma g_1^{\alpha_1}\dots
g_n^{\alpha_n}$ or as $g_1^{\alpha_1}\dots g_n^{\alpha_n}$, for
$\alpha_1,\dots,\alpha_n\in\dz$. The action of the element
$g_1^{\alpha_1}\dots g_n^{\alpha_n}$ is fixed-point-free, whereas the
action of the element $\gamma g_1^{\alpha_1}\dots g_n^{\alpha_n}$ has
a~unique fixed point, whose coordinates are
$(-\alpha_1,\dots,-\alpha_n)$.
 
We see that the action of $G_n$ on $\rr^n\sm\dz^n$ is free, and that
accordingly $(\rr^n\sm\dz^n)/G_n$ is an~$n$-dimensional manifold.  For
the points in $\dz^n\subset\rr^n$ we see that the stabilizers are of
cardinality $2$, and that given such a~point $v$, the action of its
stabilizer on an~$(n-1)$-dimensional sphere centered at $v$ is the
standard antipodal action, hence the quotient of this sphere by that
action is the projective space $\rp^{n-1}$. 

We conclude that $X_{n+1}$ can be obtained from an $n$-dimensional
manifold, whose boundary consists of $2^n$ $(n-1)$-dimensional
projective space $\rp^{n-1}$, by attaching $2^n$ cones - one at each
boundary projective space.

\subsection{The space $X_n$ as a homotopy colimit} $\,$

\nin We would like to present the space $X_n$ yet in another way, the one
which will also introduce the terminology for dealing with the case of
higher genus. The following concept is an important construction in
Combinatorial Algebraic Topology, see \cite[Chapter 15]{CAT} for
relevant background.

\begin{df} \label{hocolimdf}
The {\bf homotopy colimit}, denoted $\hocolimm\cd$, of a~diagram $\cd$
of topological spaces over a~triangulated space $\da$, is the quotient
space
\[\hocolimm\cd=\coprod_{\sigma=v_0\ra\dots\ra v_n} 
(\sigma\times\cd(v_0))/\thicksim,\]
where the disjoint union is taken over all simplices in $\da$. The
equivalence relation $\thicksim$ is generated by: for
$\tau\in\bo\sigma$, $\tau=v_0\ra\dots\ra\hat v_i\ra\dots\ra v_n$, let
$i:\tau\hra\sigma$ be the inclusion map, then
\begin{itemize}
\item for $i>0$, $\tau\times\cd(v_0)$ is identified with the subset
of $\sigma\times\cd(v_0)$, by the map induced by~$i$;
\item for $\tau=v_1\ra\dots\ra v_n$, we have 
$i(\alpha)\times x\sim\alpha\times\cd(v_0\ra v_1)(x)$, for any
$\alpha\in\tau$, and $x\in\cd(v_0)$.
\end{itemize}

\end{df}

Examples of homotopy colimits include the mapping cone and the mapping
cylinder.

For a~diagram $\cd$ over $\da$ we have a~{\it base projection map}:
\[\proj=p_b:\hocolim\cd\lra\da,\]
induced by the projections onto the first coordinate
$\sigma\times\cd(v_0)\ra\sigma$, for every simplex
$\sigma=v_0\ra\dots\ra v_n$.

\begin{figure}[hbt]
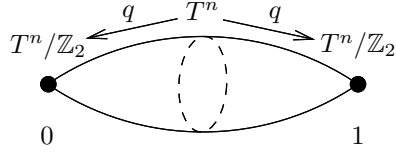

\begin{center}
  \begin{picture}(0,0)%
    \includegraphics{hcol.pstex}%
  \end{picture}%
  \input{hcol.pstex_t}%
  
\end{center}
\caption{The space $X_{n+2}$ presented as a homotopy colimit.}
\label{fig:hc1}
\end{figure}

 For $n\geq 1$, the space $X_{n+2}$ can be considered as 
a~homotopy colimit of a~diagram $\cd$ over an interval $[0,1]$, where
the latter is viewed as a~triangulated space with two vertices and one
edge. Namely, set $\cd(0):=\cd(1):=X_{n+1}=T^n/\zz$, $\cd((0,1)):=T^n$, and
let both diagram maps $\cd((0,1)\ra 0)$ and $\cd((0,1)\ra 1)$ be the
quotient maps $q:T^n\ra T^n/\zz$. See Figure~\ref{fig:hc1}.

We see that the base projection map $X_{n+2}\ra[0,1]$ is induced by 
the projection of the first coordinate of $T^{n+1}$ to the real axis:
\[(z_1,\dots,z_{n+1})\mapsto\Ree z_1.\]
By definition, this homotopy colimit is homeomorphic to the space
obtained by taking the cylinder with the base $T^n$ and then
quotioning both ends using the map~$q$. The reader is welcome to
compare this to the definition of Whitehead group, see~\cite{Co73}.

\subsection{The homology groups of $X_n$ with coefficients in $\zz$} $\,$

\nin It follows from the presentation of $X_{n+2}$ as a~homotopy colimit
that it can be seen as a~union of two disjoint copies of a~mapping
cylinder of the map $q$, glued together along the topmost copies
of~$T^n$. 

Let us slightly modify this picture. Let $b:X_{n+2}\ra[0,1]$ be the
base projection map corresponding to this diagram. We set
$A:=b^{-1}([0,1/2+\varepsilon))$ and $B:=b^{-1}((1/2-\varepsilon,1])$,
where $\varepsilon$ is a small positive number, say
$\varepsilon=0.1$. Then we have $A\cap
B=b^{-1}((1/2-\varepsilon,1/2+\varepsilon))$. Since $A\cup B=X_{n+2}$,
and the subspaces $A$, $B$, and $A\cap B$ are open in $X_{n+2}$,
we can use the Mayer-Vietoris sequence to compute the homology groups 
of $X_{n+2}$:
\begin{equation} \label{eq:mv0}
\dots\ra \wti H_i(A\cap B)\ra \wti H_i(A)\oplus \wti H_i(B)\ra 
\wti H_i(X_{n+2})\ra \wti H_{i-1}(A\cap B)\ra\dots.
\end{equation}

Clearly, both $A$ and $B$ are homotopy equivalent to $T^n/\zz$,
whereas $A\cap B$ is homotopy equivalent to $T^n$. The inclusion maps
$A\cap B\hra A$, and $A\cap B\hra B$ induce the same maps on the
homology as the quotient map $q:T^n\ra T^n/\zz$. It follows that the
Mayer-Vietoris sequence~\eqref{eq:mv0} translates to the~long exact
sequence
\begin{equation}\label{eq:mv1}
\dots\ra \wti H_i(T^n)\stackrel{(q_*,q_*)}{\lra}\wti H_i(T^n/\zz)\oplus 
\wti H_i(T^n/\zz)\ra\wti H_i(X_{n+2})\ra\wti H_{i-1}(T^n)
\stackrel{(q_*,q_*)}{\lra}\dots.
\end{equation}

This long exact sequence splits into short exact sequences as follows
from the next proposition.

\begin{prop} \label{prop:0map}
For arbitrary $n\geq i\geq 1$, the induced map 
$q_*:\wti H_i(T^n;\zz)\ra \wti H_i(T^n/\zz;\zz)$ is a~$0$-map.
\end{prop}
\pr
We know that $\wti\beta_i(T^n;\zz)={{n}\choose{i}}$. Choose
a~subset $S$ of $[n]$, such that $|S|=i$. Let $\sigma_S$ to be the sum
of all $i$-cells, indexed by $n$-tuples $(a_1,\dots,a_n)$, such that
\[
a_j=\begin{cases}
+, \text{ for } j\notin S,\\
i+ \text{ or } i-, \text{ for } j\in S,
\end{cases}
\]
for all $j=1,\dots,n$. There are $2^i$ cells like that, and their sum
is a~cycle. If we let $S$ run over all cardinality $i$ subsets of $[n]$,
we shall obtain precisely the representatives of $n\choose i$ generators
of the group $\wti H_i(T^n;\zz)$. 

All the cycles $\sigma_S$ are $\zz$-invariant, with all the $i$-cells
in each sum coming in $\zz$-invariant pairs. In every pair, both cells
are mapped to the same cell in $T^n/\zz$, hence their contributions
under the map $q_*$ cancel out each other, when we work with
$\zz$-coefficients. It follows that $q_*(\sigma_S)=0$ already on the
level of chains. Since $\sigma_s$'s generate $\wti H_i(T^n;\zz)$, we
conclude that $q_*$ is a~$0$-map on the homology groups.  \qed

\vskip5pt

Certainly, if $q_*$ is a~$0$-map, then so is the map $(q_*,q_*)$ in
the long exact sequence~\eqref{eq:mv1}. Hence, the latter splits into
the short ones of the type
\[0\ra\wti H_i(T^n/\zz;\zz)\oplus\wti H_i(T^n/\zz;\zz)\ra
\wti H_i(X_{n+2};\zz)\ra\wti H_{i-1}(T^n;\zz)\ra 0.\]
Since we are working over field coefficients we conclude that we have an
isomorphism
\[\wti H_i(X_{n+2};\zz)\cong\wti H_i(T^n/\zz;\zz)\oplus
\wti H_i(T^n/\zz;\zz)\oplus\wti H_{i-1}(T^n;\zz),\]
for all $i\geq 1$ and $n\geq 1$. Accordingly, we get
\begin{equation} \label{eq:recbet}
\wti\beta_i(T^{n+1}/\zz;\zz)=2\wti\beta_i(T^n/\zz;\zz)+
\wti\beta_{i-1}(T^n;\zz).
\end{equation}
Using that formula we arrive at the following statement.

\begin{thm} \label{thm:xnb}
The Betti numbers over $\zz$ of the space $X_{n+2}$ are given by the
following formula:
\begin{equation}\label{eq:xnb}
\wti\beta_i(X_{n+2};\zz)= 
\begin{cases}
\sum_{j=0}^{n-i+1}2^j{n-j\choose i-1},& \text{ if }n+1\geq i\geq 2;\\
0,& \text{ otherwise}.
\end{cases}
\end{equation}
\end{thm}

\pr 
Since $\wti\beta_0(T^n;\zz)=0$, the equation~\eqref{eq:recbet} implies
$\wti\beta_1(T^{n+1}/\zz;\zz)=2\wti\beta_1(T^n/\zz;\zz)$, for $n\geq
1$.  Hence
$\wti\beta_1(T^{n+1}/\zz;\zz)=2^n\wti\beta_1(T^1/\zz;\zz)=0$.

Assuming now that $i\geq 2$ we can unfold the recursive  
formula~\eqref{eq:recbet} as follows:
\begin{multline*}
\wti\beta_i(T^{n+1}/\zz;\zz)=2\wti\beta_i(T^n/\zz;\zz)+{n\choose i-1}=\\
=2\left(2\wti\beta_i(T^{n-1}/\zz;\zz)+{n-1\choose i-1}\right)+{n\choose
i-1}=\\ =4\wti\beta_i(T^{n-1}/\zz;\zz)+2{n-1\choose i-1}+{n\choose
i-1}=\\=\dots ={n\choose i-1}+2{n-1\choose i-1}+4{n-2\choose
i-1}+\dots+2^{n-i+1}{i-1\choose i-1}.\qed
\end{multline*}

For the integer coefficients the situation is somewhat more
complicated, as the map $q_*$ is not necessarilly a~$0$-map
anymore. Based on the examples for the small values of $n$ and the
general intuition for what this map is, we make the following
conjecture.

\begin{conj}
For $n\geq 1$ we have
 
\begin{equation}
\wti H_{2i}(T^n/\zz;\dz)=\zz^{a(i,n)}\oplus\dz^{b(i,n)},
\text{ for } n\geq 2i\geq 2,
\end{equation}
where $a(i,n)=\wti\beta_{2i+1}(T^n/\zz;\zz)$, and
$a(i,n)+b(i,n)=\wti\beta_{2i}(T^n/\zz;\zz)$; 
whereas $\wti H_j(T^n/\zz;\dz)=0$ for all other values of~$j$.
\end{conj}

\vskip5pt

\nin {\bf Acknowledgments.} The author would like to thank Bernd 
Sturmfels for introducing him into the world of tropical geometry and
Grisha Mikhalkin for explaining the importance of the tropical moduli
spaces. He would also like to thank Alek Vainshtein and Alex Suciu for
discussions of torus quotients, as well as Eva-Maria Feichtner for
debating various formal approaches to the topology of the space of
metric graphs. Finally, the author expresses his gratitude to the
Mathematical Institute at Oberwolfach for the hospitality during the
time a~part of this research was done.

\end{document}